\documentclass[fleqn]{article}

\usepackage{latexsym}
\usepackage{amssymb}

\usepackage{widetex2}
\usepackage{point}
\newcommand{\PTH}[1]{{\bf Proof of Theorem~#1.~}}                 
\newcommand{\PPR}[1]{{\bf Proof of Proposition~#1.~}}             
\newcommand{\PLE}[1]{{\bf Proof of Lemma~#1.}}                    
\newcommand\prothe[2]{\par\noindent\PTH{#1}{\rm #2}\medskip\par}  
\newcommand\propro[2]{\par\noindent\PPR{#1}{\rm #2}\medskip\par}  
\newcommand\prolem[2]{\par\noindent\PLE{#1}{\rm #2}\medskip\par}  

\newtheorem{thm}{Theorem}{\bfseries}{\itshape}       
\newtheorem{prop}{Proposition}{\bfseries}{\itshape}  
\newtheorem{lem}{Lemma}{\bfseries}{\itshape}         
\newtheorem{defn}{Definition}{\bfseries}{\rm}        
\newtheorem{rem}{Remark}{\itshape}{\rm}              
\newtheorem{cor}{Corollary}{\bfseries}{\itshape}     
\newtheorem{obs}{Observation}{\bfseries}{\rm}        

\def\qed{\relax\ifmmode\hskip2em \Box\else\unskip\nobreak\hskip1em $\Box$\fi}
\def\l{\ell}                                         
\def\aa{\alpha}                                      
\def\a{\mathop{\alpha}\nolimits}
\def\b{\mathop{\beta}\nolimits}
\def\G{\Gamma}                                       
\def\la{\lambda}                                     
\def\si{\sigma}                                      
\def\f{\varphi}                                      
\def\e{w}                                            
\def\FF{\mathop{\mathcal F}\nolimits}                
\def\GG{\mathop{\mathcal G}\nolimits}                
\def\TT{\mathop{\mathcal T}\nolimits}                
\def\NN{\mathop{\mathcal N}\nolimits}                
\def\RR{\mathop{\mathcal R}\nolimits}                
\def\II{\mathop{\mathcal I}\nolimits}                
\def\JJ{\mathop{\mathcal J}\nolimits}                
\def\suml {\mathop{\sum}   \limits}                  
\def\liml {\mathop{\lim}   \limits}                  
\def\capl {\mathop{\cap}   \limits}                  
\def\cdc{,\ldots,}                                   
\def\1n{1,\ldots,n}                                  
\def\0n{0,\ldots,n}                                  
\def\_#1{\mathop{\hspace{-0.45ex}^{}_{#1}}}          
\def\R{{\mathbb R}}                                  
\def\C{{\mathbb C}}                                  
\def\vj {\mathop{\tilde{J}}\nolimits}                
\def\q  {\tilde{J}}                                  
\def\J  {\tilde{J}{}}                                
\def\interca{\mathop{\scriptscriptstyle T}\nolimits} 
\def\Lpr{L^{\scriptscriptstyle\#}}                   
\def\NL{{\mathcal N}({L})}                           
\def\NJ{{\mathcal N}({\vj})}                         
\def\RL{{\mathcal R}({L})}                           
\def\RLT{{\mathcal R}({L}^{*})}                      
\def\RJ{{\mathcal R}({\vj})}                         
\def\RJT{{\mathcal R}(\J^*)}                         
\def\gri{{\scriptscriptstyle\#}}                     
\def\tr{\mathop{\rm tr}\nolimits}                    
\def\inv{\mathop{\rm inv}\nolimits}                  
\def\adj{\mathop{\rm adj}\nolimits}                  
\def\rank{\mathop{\rm rank}\nolimits}                
\def\ind{\mathop{\rm ind}\nolimits}                  
\def\beq{\begin{equation}}                           
\def\eeq{\end{equation}}                             
\def\beqq{\begin{eqnarray}}                          
\def\eeqq{\end{eqnarray}}                            
\newcommand{\abs}[1]{\left|#1\right|}                
\def\di{d}                                           
\def\ve{v}                                           
\def\tor{\to\hspace{-.07em}{*}\hspace{-.05em}}       
\def\rto{\hspace{-.02em}{*}\hspace{-.08em}\to}       
\def\ms{\mathstrut}                                  
\def\Pin{P^{\infty}}                                 

\def\baselinestretch{1.19}
\begin{document}

\title{Forest matrices around the Laplacian matrix
}

\author{Pavel Chebotarev$^1$
\ and Rafig Agaev
\\
{\normalsize Trapeznikov Institute of Control Sciences
of the Russian Academy of Sciences}
\\
{\normalsize 65 Profsoyuznaya str., Moscow 117997, Russia}
}
\footnotetext[1]{Corresponding author.\\
\indent{\ \ \it E-mail addresses:} chv@lpi.ru; pchv@rambler.ru}

\date{}

\maketitle

\begin{abstract}
We study the matrices $Q\_k$ of in-forests of a weighted digraph
$\G$ and their connections with the Laplacian matrix~$L$ of~$\G$.
The $(i,j)$ entry of $Q\_k$ is the total weight of spanning converging
forests ({\em in-forests}) with $k$ arcs such that $i$ belongs to a
tree rooted at~$j$.
The forest matrices, $Q\_k,$ can be calculated recursively and
expressed by polynomials in the Laplacian matrix; they provide
representations for the generalized inverses, the powers, and some
eigenvectors of~$L$. The normalized in-forest matrices are row
stochastic; the normalized matrix of maximum in-forests is the
eigenprojection of the Laplacian matrix, which provides an immediate
proof of the Markov chain tree theorem. A source of these results is
the fact that matrices $Q\_k$ are the matrix coefficients in the
polynomial expansion of $\adj(\la I+L)$. Thereby they are precisely
Faddeev's matrices for~$-L$.
\medskip

\noindent{\em AMS classification:\/} 05C50; 15A48
\medskip

\noindent{\em Keywords:} Weighted digraph; Laplacian matrix; Spanning
forest; Matrix-forest theorem; Leverrier-Faddeev method; Markov chain
tree theorem; Eigenprojection; Generalized inverse
\end{abstract}

\section{Introduction}
\label{sec1}

According to the matrix-tree theorem, the $(i,j)$ cofactor of the
Laplacian matrix of a weighted digraph equals the total weight of
spanning converging trees rooted at vertex $i$ of the digraph.

Fiedler and Sedl\'{a}\v{c}ek \cite{FiedlerSedlacek58} proved that the
principal minor of the Laplacian matrix resulting by the removal of the
rows and columns indexed by a set $\JJ$ is equal to the total weight of
in-forests with $\abs{\JJ}$ trees rooted at the vertices of~$\JJ$.

These results are generalized by the {\em all minors matrix tree
theorem\/} \cite{Chen76,Chaiken82} (see also~\cite{Moon94})
which expresses arbitrary minors of the Laplacian matrix in terms of
in-forests of the digraph.

We study the matrices, $Q\_k,$ of a digraph's in-forests: the $(i,j)$
entry of $Q\_k$ is the total weight of in-forests with $k$ arcs where
$i$ belongs to a tree converging to~$j$. In this paper, we show that
the forest matrices can be recursively calculated and represented by
simple polynomials in the Laplacian matrix~$L$; in turn, the powers of
$L$ are linear combinations of~$Q\_k$'s. Further, we demonstrate that
the forest matrices are useful to interpret a number of expressions
that involve the Laplacian matrix, including those of the group and
Moore-Penrose inverses, and some eigenvectors. Of special interest is
the normalized matrix $\J$ {\em of maximum in-forests\/} of a digraph
previously used~\cite{LeightonRivest83,LeightonRivest86} to represent
the long run transition probabilities of Markov chains. We prove that
$\J$ is the eigenprojection of the Laplacian matrix corresponding to
the eigenvalue~$0$ and study some properties of~$\J.$

A seminal result that enables one to give short algebraic proofs to
these representations is the fact that matrices $Q\_k$ coincide with
the matrix coefficients in the polynomial form of $\adj(\la I+L)$:
\[
\adj(\la I+L)=\suml_{k=0}^{n-1}Q\_{n-k-1}\la^k,
\]
where $\adj A$ is the transposed matrix of cofactors of~$A$. This
expansion is a corollary to the {\em parametric matrix-forest
theorem\/}~\cite{AgaChe00} which expresses the entries of $(I+\tau
L)^{-1},$ $\tau\in\R$ in terms of in-forests.

All results of this paper are applicable to unweighted digraphs
(by taking all weights equal to one) and undirected graphs (by
considering symmetric digraphs).

The paper is organized as follows. After the notation section, we
briefly survey the major known results on the minors of the Laplacian
(Kirchhoff) matrix of a weighted digraph (Section~\ref{Preli}), give
a new proof
to the matrix-forest
theorem for digraphs (Section~\ref{nata}), present a recursive method
for calculating the forest matrices (Section~\ref{s_calc}), establish
polynomial representations of the forest matrices
(Sections~\ref{s_rela}), study the normalized matrix $\J$ of
maximum in-forests (Section~\ref{s_maxi}), consider $L$ and $\J$ as
linear transformations and show that $\J$ is the eigenprojection of
$L,$ which yields the {\em Markov chain tree theorem\/}
(Section~\ref{line}), and finally, express the generalized inverses of
$L$ in terms of the forest matrices (Section~\ref{pseudo2}).

\section{Notation}
\label{sec2}

\subsection{Graph definitions}

For graph terminology, we mainly follow~\cite{Harary69}.
Suppose that $\G$ is a weighted digraph without loops, $V(\G)=\{\1n\},$
$n>1$, is its set of vertices and $E(\G)$ its set of arcs. The
weights of all arcs are strictly positive.
Let $W=(\e\_{ij})$ be the matrix of arc weights of~$\G$. Its $(i,j)$
entry, $\e\_{ij},$ equals zero iff there is no arc from vertex
$i$ to vertex~$j$ in~$\G$. If $\G^{\prime}$ is a subgraph of $\G$, then
the weight of $\G^{\prime}$, $\e(\G^{\prime})$, is the product of the
weights of all its arcs; if $\G'$
does not contain arcs, then $\e(\G')=1$. The weight of a nonempty set
of digraphs $\GG$ is defined as follows:
\beq
\label{set_weight}
\e(\GG)=\suml_{H\in\GG}\e(H);\quad \e(\varnothing)=0.
\eeq

A~{\it spanning\/} subgraph of $\G$ is a subgraph of $\G$ with vertex
set $V(\G)$. The {\it outdegree\/}
of vertex $\ve$ is the number of arcs that come from~$\ve$.
A {\it converging tree\/} is a weakly connected (i.e., its
corresponding undirected graph is connected) digraph in which one
vertex, called the {\it root}, has outdegree zero and the remaining
vertices have outdegree one.

A~converging tree is said to {\em converge\/} to its root. Spanning converging
trees are sometimes called {\it in-arborescences}.
A~{\it converging forest\/}
is a digraph all of whose weak components (i.e., maximal weakly
connected subgraphs) are converging trees. The roots of these trees are
the roots of the converging forest.

\begin{defn}
{\rm
An {\em in-forest\/} is a spanning converging forest.
}
\end{defn}

\begin{defn}
\label{De2}
{\rm
An in-forest $F$ of a digraph $\G$ is called a {\em maximum
in-forest\/} of $\G$ if $\G$ has no in-forest with a greater number of
arcs than in~$F$.
}
\end{defn}

{\em Out-forests\/} which {\em diverge\/} from their roots
and {\em maximum out-forests\/} are defined in the same manner. In this
paper, we deal with in-forests, but a parallel theory can be
developed for out-forests.

The notion of maximum in-forest of a digraph generalizes the
concept of spanning converging tree (in-arborescence). If spanning
converging trees of a digraph exist, they coincide with maximum in-forests;
otherwise maximum in-forests inherit some of their properties.
These properties were studied in~\cite{AgaChe00}.

It is easily seen that every maximum in-forest of $\G$ has the minimum
possible number of converging trees; we call this number the {\it
in-forest dimension\/} of $\G$ and denoted it by~$\di$. The number of
arcs in any maximum in-forest is obviously $n-\di$; in general, the
number of disjoint trees in a spanning forest with $k$ arcs is $n-k$.

By $\FF^{\tor}(\G)=\FF^{\tor}$ and $\FF^{\tor}_k(\G)=\FF^{\tor}_k$
we denote the set of all in-forests of $\G$ and the set of all
in-forests of $\G$ with $k$ arcs, respectively; $\FF^{i\tor j}_k$
will designate the set of all in-forests with $k$ arcs where $i$
belongs to a tree converging to~$j$; $\FF^{i\tor
j}=\cup_{k=0}^{n-\di}\FF^{i\tor j}_k$ is the set of such
in-forests with any number of arcs. The notation
$\FF^{\tor}_{(k)}$ will be used for the set of in-forests that
consist of $k$ trees. The $\to\hspace{-.5em}{*}$ sign relates to
in-forests; the corresponding notation for out-forests is
$\FF^{\rto},$ etc.

Let
\beqq
\label{sik}
\si\_k
&=&\e(\FF^{\tor}_k),\quad k=0,1,\ldots,
\\
\si
&=&\e(\FF^{\tor})=\suml_{k=0}^{n-\di}\si\_k.
\label{si}
\eeqq

By (\ref{sik}) and (\ref{set_weight}), $\si\_k=0$ whenever $k>n-\di,$
and $\si\_0=1.$

We will also consider the parametric value
\beq
\label{sitau}
\si(\tau)
=\suml_{k=0}^{n-\di}\si\_k\tau^k,
\eeq
which is the total weight of in-forests in $\G$ provided that all arc
weights are multiplied by~$\tau$.

Let
\beq
\label{sk}
s\_k=\suml_{j=0}^k\si\_j,\quad k=\0n-\di
\eeq
be the total weight of in-forests of $\G$ with at most $k$ arcs. Then,
by definition, $s\_{n-\di}=\si$.

Finally,
\beq
\label{sktau}
s\_k(\tau)=\suml_{j=0}^k\si\_j\tau^j,\quad k=\0n-\di,
\eeq
whence $s\_{n-\di}(\tau)=\si(\tau).$

\subsection{Matrix definitions}

For any $n\!\times\!n$ matrix $A$, let
$A(\overline{\II}\,|\,\overline{\JJ}),$ where $\II,\JJ\subseteq \{1\cdc
n\},$ be the submatrix of $A$ obtained by the removal of the rows
indexed by $\II$ and the columns indexed by~$\JJ$. For a complex matrix
$A,$ $A^*$ is the conjugate transpose (Hermitian adjoint) and
$A^{\interca}$ the transpose of~$A.$

The {\it Laplacian\/} (or {\it row Laplacian}) matrix of a
weighted digraph $\G$ is the $n\!\times\!n$ matrix $L=L(\G)=(\l\_{ij})$
with entries $\l\_{ij}=-\e\_{ij}$ when $j\ne i$ and
$\l\_{ii}=-\suml_{k\ne i}\l\_{ik}$, $i,j=\1n$. The {\it column
Laplacian} matrix $L'=L'(\G)=(\l'_{ij})$ differs from $L$ by the
diagonal only: ${\l'_{ij}=-\e\_{ij}}$ when $j\ne i$ and
$\l'_{ii}=-\suml_{k\ne i}\l'_{ki}$, $i,j=\1n$. The {\it Kirchhoff\/}
(or {\it row Kirchhoff\/}) matrix \cite{Tutte84} is $K=L'^{\interca};$
the {\it column Kirchhoff} matrix is $K'=L^{\interca}.$ These four
singular matrices are generalizations of the Laplacian (Kirchhoff)
matrix of an undirected graph. In what follows, we deal with the
Laplacian matrix~$L(\G)$ and reformulate for it some results originally
obtained for the other matrices.

Throughout let $\G$ be a fixed digraph.
Consider the {\em matrices
\[
Q\_k=(q_{ij}^k),\quad k=0,1,\ldots,
\]
of in-forests of\/ $\G$ with $k$ arcs}: the entries of $Q\_k$ are
\beq
\label{qijk}
q_{ij}^k=\e(\FF_k^{i\tor j}).
\eeq

By (\ref{qijk}) and (\ref{set_weight}), $Q\_k=0$ whenever $k>n-\di,$ and
$Q\_0=I.$

The {\em matrix of all in-forests\/} is
\beq
\label{Q}
Q=(q\_{ij})=\suml_{k=0}^{n-\di}Q\_k
\eeq
with entries $q\_{ij}=\e(\FF^{i\tor j}).$

We will also consider the {\em normalized matrices of forests}:
\beq
\label{Jk}
J\_k=\si_k^{-1}Q\_k,\quad k=\0n-\di,
\eeq
\beq
\label{J}
J=\si^{-1}Q,
\eeq
and the parametric matrices
\beq
\label{Qtau}
Q(\tau)=\suml_{k=0}^{n-\di} Q\_k\tau^k,
\eeq
\beq
\label{Jtau}
J(\tau)=\si^{-1}(\tau)\,Q(\tau),\quad \tau\ge0,
\eeq
where $\si\_k,$ $\si,$ and $\si(\tau)$ are defined by
(\ref{sik})--(\ref{sitau}).

The {\em normalized matrix of maximum in-forests\/} $J\_{n-\di}\/$
will be also denoted by~$\J$:
\[
\J=J\_{n-\di}.
\]
In the case of undirected graphs, the entries of $\J$ are the same
within every connected component. In the directed case, this
matrix possesses nontrivial properties determined by the properties of
maximum in-forests, cf.~\cite{AgaChe00}.

\begin{prop}
\label{pro_sto}
The matrices $J\_k,\;k=\0n-\di,$ $J,$ and $J(\tau)$ are row stochastic.
\end{prop}

\propro{\ref{pro_sto}}{
Every row sum of $Q\_k,\;k=\0n-\di,$ is $\si\_k$. Indeed, for every
$i=\1n,$ we have
\[
\suml_{j=1}^nq^k_{ij}
=\suml_{j=1}^n\e(\FF_k^{i\tor j})
\stackrel{\scriptstyle{(\!\ast\!)}}{=}
 \e\Bigl(\bigcup_{j=1}^n\FF_k^{i\tor j}\Bigr)
=\e(\FF^{\tor}_k)
=\si\_k.
\]
In the $(\ast)$ passage, we used the fact that
$\FF_k^{i\tor j\_1}\cap\FF_k^{i\tor j\_2}=\varnothing$ whenever
$j\_1\ne j\_2.$ Thus, the nonnegative matrices $J\_k=\si_k^{-1}Q\_k$
are row stochastic. Now the stochasticity of $J$ and $J(\tau)$ follows from
their definitions.
\qed
}

The aim of this paper is to interpret, in terms of the forest matrices,
a number of expressions that involve the Laplacian matrix
as well as to provide polynomial expressions for the forest matrices
themselves.

\section{Preliminaries}
\label{Preli}

This section briefly surveys some known results on the minors of
a digraph's Laplacian matrix.

The oldest result of this kind is the {\em matrix-tree theorem\/}
by Tutte~\cite{Tutte48,Tutte84}, although some authors (e.g.,
\cite{Chaiken82}; cf.\ \cite{Moon70}) trace it back to
Sylvester~\cite{Sylvester57} and its proof to Borchardt~\cite{Borchardt60}.

\begin{thm}
\label{MTT}
For every $i,j\in V(\G),$ $\ell^{ij}=\e(\TT^{\tor i})$ holds$,$ where
$\ell^{ij}$ is the cofactor of the $(i,j)$ entry of $L$ and $\TT^{\tor
i}$ is the set of all spanning trees converging to $i$ in~$\G$.
\end{thm}

As stated in \cite{KelmansPakPostnikov99}, ``This small formula opens
a world of opportunities.''

Tutte~\cite{Tutte84} formulated this theorem for the diagonal cofactors
of the Kirchhoff matrix. A~version that involves all
cofactors of the Laplacian and the column Laplacian matrices can be
found in~\cite{Harary69}.
We do not describe multiple analogues of the matrix-tree theorem here.

By definition, $L$ has the form $L=D-W,$ where $W$ is the nonnegative
matrix of arc weights and $D$ is the diagonal matrix ensuring the zero
row sums of~$L.$ Therefore, by Ger\v{s}gorin's theorem, the real part
of each nonzero eigenvalue of $L$ is positive. Thus, $L$ is a singular
M-matrix (see, e.g., \cite[Theorem~4.6 in Chapter~6]{BermanPlemmons79}).
One of the consequences is that all the principal minors of $L$ are
nonnegative. Fiedler and Sedl\'{a}\v{c}ek~\cite{FiedlerSedlacek58}
obtained an interpretation of all principal minors of the Laplacian
matrix in terms of spanning forests:

\begin{thm}
\label{FiSe}
For any $\JJ\subseteq\{\1n\}$, $\det
L(\overline{\JJ}\,|\,\overline{\JJ})=\e(\FF^{\tor \JJ})$ holds$,$
where $\FF^{\tor \JJ}$ is the set of in-forests for which $\JJ$ is
the set of roots.
\end{thm}

Later this theorem was formulated and proved in~\cite{ChaikenKleitman78}.
Its special case with undirected graphs and $\abs{\JJ}=2$ was discovered
and employed earlier in the theory of electrical networks (see,
e.g.,~\cite{Percival53}). Fiedler and Sedl\'{a}\v{c}ek stated their
result for the column Laplacian matrix and out-forests. Generally, to
get interpretations for the minors of the column Laplacian matrix
$L'(\G),$ it suffices to observe that for the digraph obtained from
$\G$ by the reversal of all arcs, the Laplacian matrix coincides with
$K(\G)=L'^{\interca}(\G)$ and the in-forests are in a weight preserving
correspondence with the out-forests of~$\G.$

Let
\beq
\label{p(la)}
\f(\la)=\det(\la I+L)=\suml^n_{k=0}c\_{n-k}\la^k
\eeq
be the characteristic polynomial of~$-L$ and let $\si\_k$ be as
defined in~(\ref{sik}).

\begin{prop}
\label{char-L}
In\/ $(\ref{p(la)}),\;$ $c\_k=\si\_k,\;$ $k=0,\ldots,n.$
\end{prop}

In view of Theorem~\ref{FiSe}, this proposition follows from
the fact that $c\_k$ is equal to the sum of the $k\!\times\!k$
principal minors of~$L$. In the case of undirected unweighted
multigraphs, Proposition~\ref{char-L} is due to
Kelmans~\cite{Kelmans67,KelmansChelnokov74}, who was probably the first
\cite{Kelmans6566E} to study the Laplacian characteristic polynomial
(see also discussion in \cite[p.~42]{Moon70} and \cite[Sections
1.2,~1.5]{CvetkovicDoobSachs80}, and \cite[Theorem~7.5]{Biggs74}); some
extensions are given in~\cite[the last statement on
p.~236]{BapatConstantine92} and~\cite[Theorem~2]{ChungLanglands96}. An
alternative representation for the coefficients of the Laplacian
characteristic polynomial can be found in~\cite{Forman93}.

Since $\si\_k=0$ if and only if $k>n-\di$ ($k=0,1,\ldots$),
Proposition~\ref{char-L} implies

\begin{cor}
\label{mult0}
The multiplicity of\/ $0$ as the eigenvalue of $L$ is~$\di.$
\end{cor}

Another immediate consequence of Proposition~\ref{char-L} is
\smallskip

\begin{cor}
\label{sumprod}
${\displaystyle
\:\si\_k=\sum_{\JJ:\:\abs{\JJ}=k}\:\prod_{j\in\JJ}\la\_j,
\quad k=0,\ldots,n,
}$\\
where $\la\_1\cdc\la\_n$ are the eigenvalues of $L$ and $\JJ$ are
the subsets of\/ $\{1\cdc n\}$.
\end{cor}

Chen~\cite[p.~313, Problems~4.14 and 4.16]{Chen76} proposed an
extension of the matrix-tree theorem to additional minors of the
Laplacian matrix and Chaiken~\cite{Chaiken82} gave a similar graph
interpretation to all minors of~$L'.$ Moon~\cite{Moon94}
obtained a more general expansion which applies to all minors of
arbitrary matrices; Chaiken's theorem and a number of W.K.~Chen's
expansions follow from his result as special cases.
Minoux~\cite{Minoux99} generalized Chaiken's theorem to semirings and
Bapat et al.~\cite{BapatGrossman99} to mixed graphs (where each arc is
either directed or undirected). Other useful graph interpretations of
minors and determinants are given in~\cite{MaybeeOlesky89}.

We do not quote these results here, but we employ Chaiken's formulation
\cite{Chaiken82} of the {\it all minors matrix tree theorem\/} in the
proof of a {\em matrix-forest theorem\/} in the following section.

\section{Another matrix-forest theorem}
\label{nata}

The following theorem~\cite{CheSha95a,CheSha97} provides expressions
for the forest matrices $Q$ and $J$ (see (\ref{Q}) and~(\ref{J})) in
terms of the cofactors and the determinant of $I+L,$ where $I$ is the
identity matrix.

\begin{thm}
\label{mft0}
$Q=\adj(I+L)$ and $\si=\det(I+L).$
Thus$,$ $J=(I+L)^{-1}.$
\end{thm}

For the properties of $(I+L)^{-1},$
see~\cite{CheSha97,CheSha98,Merris97,Merris98}.

It was mentioned in \cite{CheSha97} that a quick way to prove the
matrix-forest theorem is to employ the all minors matrix tree theorem,
more specifically, to apply the first formula (without number) on
page~328 in~\cite{Chaiken82}. Below we give a complete inference of
Theorem~\ref{mft0} from the all minors matrix tree theorem.
Note that a self-contained
proof of the matrix-forest theorem for unweighted
multigraphs can be found in~\cite{Sha94}. Another inference based on
some results of
\cite{Kelmans67,KelmansChelnokov74,FiedlerSedlacek58,MaybeeOlesky89}
was given in~\cite{CheSha95a} for the case of weighted multidigraphs
and multigraphs. Undirected and unweighted analogies of
Theorem~\ref{mft0} have been presented in~\cite{Merris97,Merris98}
(with the proof based on Chaiken's theorem) and~\cite{CheSha95}.

In the following proof of Theorem~\ref{mft0}, we employ a standard
trick which enables one to reduce many novel statements about forests
to known statements about trees or forests.
Versions of this trick have been used in many papers, e.g.,
\cite{BapatConstantine92,Chaiken82,CheSha01,Chen76,GolenderDrboglav81,%
JonesPittel99,KelmansPakPostnikov99,KirklandNeumann99b,Merris97,Merris98,%
ProppWilson98}.%
\footnote{Note that one more expedient is to identify the roots of all
trees in a forest, which converts the forest into a tree
\cite{CvetkovicDoobSachs80,FiedlerSedlacek58,Kelmans67,KelmansChelnokov74,%
Myrvold92,CheSha95a}.}
We formalize it by
\begin{defn}
{\rm
Let $\G$ be a weighted digraph. The digraph $\widehat{\G}$ with vertex
set $V(\widehat{\G})=V(\G)\cup\{0\},$ arc set
$E(\widehat{\G})=E(\G)\cup\{(j,0)\,:\,j\in V(\G)\},$ the weights of
arcs in $E(\widehat{\G})\cap E(\G)$
the same as for~$\G,$ and $\e((j,0))=1,\;j\in V(\G),$ will be called
the {\em ground extension\/} of~$\G$.\footnote{In
\cite{KelmansPakPostnikov99} $\widehat{\G}$ is called the {\em cone\/}
of~$\G$.}
}
\end{defn}

\begin{obs}
\label{obs1}
Let $\widehat{\G}$ be the ground extension of\/ $\G.$ Let $U=I+L(\G),$
$\widehat{L}=L(\widehat{\G}).$ Then for any $\II,\JJ\subseteq V(\G),$
$U(\overline{\II}\,|\,\overline{\JJ})
=\widehat{L}(\overline{\II\cup\{0\}}\,|\,\overline{\JJ\cup\{0\}})$ holds.
\end{obs}

By virtue of Observation~\ref{obs1}, if one has expressions for all
minors of the Laplacian matrices $L$ (say, those provided by the
all minors matrix tree theorem), then expressions for all minors of
matrices $I+L$ are got gratis. The following lemma establishes a
correspondence between the forests in $\G$ and some forests in~$\widehat{\G}$.
The lemma is formulated here in a form useful for expressing all minors of
$I+L.$

\begin{lem}
\label{trick}
Consider
$\II=\{i_1\cdc i_k\}\subseteq V(\G),$
$\JJ=\{j_1\cdc j_k\}\subseteq V(\G),$
$0\le k\le n,$ and the set of in-forests
$\FF^{\tor}\cap\Bigl(\capl_{u=1}^k\FF^{i\_u\tor j\_u}\Bigr)$ in $\G$.
Then there exists a weight preserving one-to-one correspondence between
this set and the set
$\widehat{\FF}\ms^{\tor}_{\II\!\JJ}$ of in-forests
$F\in\widehat{\FF}\ms^{0\tor 0}\cap\Bigl(\capl_{u=1}^k\widehat{\FF}\ms^{i\_u\tor j\_u}\Bigr)$
in $\widehat{\G}$ such that the $F$'s consist of exactly $k+1$ trees.
\end{lem}

\prolem{\ref{trick}}{
Let $F\in\FF^{\tor}\cap\Bigl(\capl_{u=1}^k\FF^{i\_u\tor j\_u}\Bigr)$.
To define the corresponding forest in
$\widehat{\FF}\ms^{\tor}_{\II\!\JJ},$ consider the replica $F'$ of $F$
in $\widehat{\G}$ and attach the arcs $(r,0)$ to it, where the $r$'s
are the roots of $F'$ that are not in $\JJ$. The resulting in-forest
consists of exactly $k+1$ trees and belongs to
$\widehat{\FF}\ms^{\tor}_{\II\!\JJ}$. Conversely, for any
$\widehat{F}\in\widehat{\FF}\ms^{\tor}_{\II\!\JJ},$ consider its
restriction to $V(\G)$ as the corresponding forest of $\G$. Obviously,
this correspondence is one-to-one and the corresponding forests share
the weight.
\qed
}

\prothe{\ref{mft0}}{
Consider the ground extension $\widehat{\G}$ of $\G$.
By Observation~\ref{obs1}, if $U=I+L(\G),$
$U^{ij}$ is the $(i,j)$ entry of~$\adj U,$ and $\widehat{L}=L(\widehat{\G}),$
then
\beq
\label{26}
U^{ij}
=(-1)^{i+j}\det U(\overline{\{j\}}\,|\,\overline{\{i\}})
=(-1)^{i+j}\det \widehat{L}(\overline{\{0,j\}}\,|\,\overline{\{0,i\}}).
\eeq

Let $\widehat{\FF}\ms^{0\tor 0,i\tor j}_{(2)}$ be the set of in-forests
$F\in\widehat{\FF}\ms^{0\tor 0}\cap\widehat{\FF}\ms^{i\tor j}$ that
consist of two trees. Denoting by $\inv\{0\to 0,i\to j\}$ the number
of violations of monotonicity in the two-element correspondence $\{0\to
0,i\to j\}$, which is obviously zero, and using the all minors matrix
tree theorem~\cite{Chaiken82,Moon94}, we get
\beqq
\!\!\!\!\!\!
\det\widehat{L}(\overline{\{0,j\}}\,|\,\overline{\{0,i\}})
&\!\!\!\!=\!\!\!\!&
     (-1)^{\abs{\{k\in V(\G)\,:\,k<j\}}\,+\,\abs{\{k\in V(\G)\,:\,k<i\}}}
     \!\!\!\!\!\!\!\!\!\!\!\!
     \suml_{F\in\widehat{\FF}\ms^{0\tor 0,i\tor j}_{(2)}}
     \!\!\!\!\!\!\!\!\!\!
     (-1)^{\inv\{0\to 0,i\to j\}}\e(F)
\cr
\label{27}
&\!\!\!\!=\!\!\!\!&(-1)^{j+i-2}\e(\widehat{\FF}\ms^{0\tor 0,i\tor j}_{(2)}).
\eeqq
In the first passage, we used the fact that
$\widehat{\FF}\ms^{i\tor 0}=\varnothing.$

Lemma~\ref{trick} implies $\e(\widehat{\FF}\ms^{0\tor 0,i\tor
j}_{(2)})=\e(\FF^{i\tor j}),$ so, from (\ref{26}) and (\ref{27}), we get
\[
U^{ij}=(-1)^{2i+2j-2}\e(\FF^{i\tor j})=\e(\FF^{i\tor j})=q\_{ij}.
\]

By Observation~\ref{obs1}, Theorem~\ref{MTT}, and Lemma~\ref{trick},
$
\det U
=\det\widehat{L}(\overline{\{0\}}\,|\,\overline{\{0\}})
=\e(\widehat{\FF}\ms^{0\tor 0}_{(1)})
=\e(\FF^{\tor})
=\si
$
(cf.\ \cite[Eq.~(37)]{JonesPittel99} and
      \cite[7.2 and~7.3]{KelmansPakPostnikov99}).
This completes the proof.
\qed
}

\begin{rem}
\label{rem_pos}
{\rm
Obviously, the positivity of arc weights is needed for the last
statement of Theorem~\ref{mft0} only; the first two statements are
preserved for digraphs with arbitrary arc weights.
}
\end{rem}

\begin{rem}
\label{r_drbog}
{\rm
Note that the cofactors and the determinant of $I+L,$ in the case of an
unweighted undirected graph $G$, have been expressed in
\cite{GolenderDrboglav81} in terms of spanning trees and 2-forests in
the ground extension of $G$ (for the case of weighted graphs,
cf.~\cite[Theorem~2.3]{KirklandNeumann97}).
Ref.~\cite{GolenderDrboglav81} also discusses the idea of using graph
invariants related to $(I+L)^{-1}$ in the study of the graph
isomorphism problem. We surmise that the forest matrices $Q\_k$ also
have some potential in this respect.
}
\end{rem}

It is easily seen that $I+\tau L$ with $\tau\ge0$ are nonsingular
M-matrices, so their inverses are nonnegative.
In the next section, the following parametric matrix-forest
theorem~\cite{AgaChe00} will be helpful:
\bigskip

\noindent
{\bf Theorem~\ref{mft0}$'.$\ }
{\it For
any $\tau\in\R,$
$Q(\tau)=\adj(I+\tau L)$ and
$\si(\tau)=\det(I+\tau L).$ Thus$,$ for any $\tau\ge0,$
$J(\tau)=(I+\tau L)^{-1}.$
}
\bigskip

To prove this theorem, it suffices to apply Theorem~\ref{mft0}
to the weighted digraph $\G'(\tau)$ that differs from
$\G$ in the weights of arcs only: for all $i,j=\1n,$
$\e'_{ij}(\tau)=\tau\e\_{ij}$.
By Remark~\ref{rem_pos}, the nonnegativity of $\tau$ is needed
for the last statement of Theorem~\ref{mft0}$'$ only.

\section{A method for calculating $Q\_1\cdc Q\_{n-\di}$}
\label{s_calc}

We first show that $Q\_1\cdc Q\_{n-\di}$ are the matrix
coefficients in the polynomial expansion of
$\adj(\la I+L)$.
\vspace{-1em}

\begin{prop}
\label{MatCo}
${\displaystyle
\adj(\la I+L)
=\suml_{k=0}^{n-d}Q\_k\la^{n-k-1}.}
$
\end{prop}

\propro{\ref{MatCo}}{
If $\la=0,$ then the right-hand side is zero whenever $\di>1$ and
it reduces to $Q\_{n-1}$ when $d=1$ (we put $\la^0\equiv1$). This is
equal to $\adj(\la I+L)$ by Theorem~\ref{MTT}. For any $\la\ne0,$ let
$\tau=\la^{-1}$. Using Theorem~\ref{mft0}$'$ we get
\beq
\label{adjexp}
\adj(\la I+L)
=\adj\la(I+\tau L)
=\la^{n-1}Q(\tau)
=\la^{n-1}\suml_{k=0}^{n-\di} Q\_k\tau^k
=\suml_{k=0}^{n-\di} Q\_k\la^{n-k-1}.
\!\!\!\!\!\!\!\!\qed
\eeq
}
\vspace{-1em}

Proposition~\ref{MatCo} underlies an easy algorithm for calculating
$Q\_1\cdc Q\_{n-\di}$ and $\si\_1\cdc\si\_{n-\di}.$

\begin{prop}
\label{pro.allk}
For any $k=0,1,\ldots,$
\beqq
\label{Fadd1}
Q_{k+1}
\!\!&=\!\!&(-L)Q_{k}+\si\_{k+1}\!I,\\
\label{Fadd2}
\si_{k+1}
\!\!&=\!\!&\frac{\tr(LQ\_k)}{k+1}.
\eeqq
\end{prop}

\propro{\ref{pro.allk}}{
Since, by Proposition~\ref{MatCo}, $Q\_0\cdc Q\_{n}$ are the matrix
coefficients in the polynomial form of $\adj(\la I+L),$ where $\la I+L$
is the characteristic matrix of $-L$ and, by Proposition~\ref{char-L},
$\si\_0\cdc \si\_{n}$ are the coefficients of the characteristic
polynomial of $-L,$ the equations~\cite[\S3 of Chapt.~4]{Gantmacher66}
$Q\_{k+1}=\si\_{k+1}\!I-LQ\_k,\;\; k=0,1,\ldots,$ take place.

To prove (\ref{Fadd2}), it suffices to take the traces on the left and
on the right of (\ref{Fadd1}) and use the fact that
\[
\tr Q_{k}=(n-k)\si\_{k},\quad k=0,1,\ldots,
\]
which holds since every in-forest with $k$ arcs has $n-k$ roots.
\qed
}

Note that, by virtue of Propositions~\ref{char-L} and~\ref{MatCo}, the
recurrent application of (\ref{Fadd2}) and (\ref{Fadd1}) starting with
$Q\_0=I$ coincides with the Leverrier-Faddeev
algorithm~\cite{FaddeevFaddeeva59,Gantmacher66} applied to calculate
the characteristic polynomial of~$-L$.

Consider now a few corollaries to Proposition~\ref{pro.allk}. First,
in what follows we will need a recurrence formula for the row stochastic
matrices $J_{k}.$ It is:
\beq
\label{j_k+1}
J_{k+1}=\frac{\si\_k}{\si\_{k+1}}(-L)J_{k}+I, \quad k=0\cdc n-\di-1.
\eeq

Second, the matrices $L Q_{k}$ prevailing in Proposition~\ref{pro.allk}
have a noteworthy graph interpretation. Let $\G\_{k}$ be the {\em
digraph of in-forests with $k$ arcs\/} of $\G$, i.e., the digraph on
vertex set $V(\G\_{k})=V(\G)$ whose matrix of arc weights results from
$Q_{k}$ by putting zeros on the main diagonal. In other words,
$(i,j)\in E(\G\_k)$ whenever $j\ne i$ and $q_{ij}^k>0$; $q_{ij}^k$ is
the weight of such arc. Evidently, $\G\_1=\G.$

\begin{prop}
\label{allmatrbek}
$L Q_{k}$ is the Laplacian matrix of\/~$\G\_{k+1},$
$k=0,1,\ldots.$
\end{prop}

\propro{\ref{allmatrbek}}{ By Proposition~\ref{pro.allk},
$LQ_{k}=\si\_{k+1}\!I-Q\_{k+1},$ so the off-diagonal entries of
$LQ\_k$ coincide with those of $L(\G\_{k+1}).$ To complete the proof,
note that every row sum of $LQ\_k$ is zero, since every row sum of both
$\si\_{k+1}\!I$ and $Q\_{k+1}$ is~$\si\_{k+1}.$
\qed
}

Finally, Proposition~\ref{pro.allk} provides a recurrent
formula for the Laplacian matrices $L\_{k}\!:=\!L(\G\_{k})$:
\[
L\_{k+1}=L\,\Bigl(-L\_k+\frac{\tr L\_k}{k}I\Bigr),\quad k=1,2,\ldots.
\]

We are going to discuss the application of digraphs $\G\_k$ to the
analysis of $\G$ elsewhere.

\def\baselinestretch{1.08}

\section{Forest matrices as polynomials in the Laplacian\\ matrix}
\label{s_rela}

It follows from Proposition~\ref{pro.allk} that the forest matrices
$Q\_{k},$ $Q,$ and $Q(\tau)$ are polynomials in~$L.$ As a corollary,
the powers of $L$ are linear combinations of $Q\_0\cdc Q\_{n-\di}.$

First, it is straightforward to prove
\def\nesk{\vspace{-1.1ex}}
\nesk

\begin{prop}
\label{teo.allk}
${\displaystyle Q\_{k}=\suml_{i=0}^k\si\_{k-i}(-L)^i,\;\;\; k=0,1,\ldots.}$
\end{prop}
\nesk

These expressions are closely related to the characteristic polynomial
of~$-L$ (\ref{p(la)}) which, by Proposition~\ref{char-L}, can be
represented as
$
\f(\la)=(...((\si\_0\la+\si\_1)\la+\si\_2)\la+\ldots+\si\_{n-1})\la+\si\_n.
$
To find $\f(\la),$ one can successively calculate
$\f\_0(\la)=\si\_0,$
$\f\_1(\la)=\si\_0\la+\si\_1,$
$\f\_2(\la)=(\si\_0\la+\si\_1)\la+\si\_2\cdc$
$\f\_n(\la)=\f(\la).$
It is easily seen now that $Q\_k=\f\_k(-L),$ $k=0\cdc n.$
\vspace{-.7ex}

\begin{cor}
\label{cor_comm}
{The matrices $Q\_k,$ $k=0,1,\ldots,$ commute with all matrices with
which $L$ commutes$,$ in particular$,$ with $L,$ $Q(\tau),$ and each
other.
}
\end{cor}
\vspace{-.7ex}

By Theorems~\ref{mft0} and~\ref{mft0}$',$ $Q=\adj(I+L)$ and
$Q(\tau)=\adj(I+\tau L).$ Proposition~\ref{teo.allk}, (\ref{Q}),
and (\ref{Qtau}) provide a polynomial form of $Q$ and $Q(\tau).$

\label{Sec_Poli}
\begin{prop}
\label{Qexpan}
\beqq
Q
&\!=&\!\suml_{k=0}^{n-\di}s\_{n-\di-k}(-L)^k
     =\adj(I+L),
\cr
\label{iden3}
Q(\tau)
&\!=&\!\suml_{k=0}^{n-\di}s\_{n-\di-k}(\tau)\,(-\tau L)^k
     =\adj(I+\tau L),
\eeqq
where $s\_i$ and $s\_i(\tau)$ are defined in\/ $(\ref{sk})$
and\/~$(\ref{sktau}).$
\end{prop}

By (\ref{adjexp}), $\adj(\la I+L)=\la^{n-1}Q(\tau),$ where $\la\ne0$
and $\tau=1/\la.$ Combining this with (\ref{iden3}) and (\ref{sktau}),
we obtain
\nesk\vspace{-.7ex}

\begin{cor}
\label{dualadj}
${\displaystyle
\adj(\la I+L)
=\suml_{k=0}^{n-d}s'_{n-d-k}(\la)\,(-L/\la)^k,
}$\\
where
$s'_i(\la)=\sum_{j=0}^i\si\_j\la^{n-j-1},\;\;i=\0n-\di,$ and\/ $\la\ne0.$
\end{cor}
\vspace{-.9ex}

Corollary~\ref{dualadj} and Proposition~\ref{MatCo} can be considered
as dual representations of $\adj(\la I+L)$.

It follows from Proposition~\ref{teo.allk} that the powers of $L$
are linear combinations of $Q\_0\cdc Q\_{n-\di},$
but the coefficients are more complicated than before.

\begin{prop}
\label{pro-combi}
For $m=0,1,\ldots,\;\;$
${\displaystyle(-L)^m=\suml_{k=0}^m\a\_kQ\_{m-k}}$ holds$,$ where
$\a\_0=1,$
\beq
\label{e_combi}
\a\_k
=\suml_{(p\_1\cdc p\_k):\;\sum ip\_i=k}
(-1)^{\sum p\_i}
\frac{\bigl(\sum p\_i\bigr)!}{\prod\bigl(p\_i!\bigr)}
\prod\si_i^{p\_i},\quad k=1\cdc m,
\eeq
$p\_i$ are nonnegative integers$,$ and all sums and products in
$(\ref{e_combi}),$ except for the first sum$,$ range from $i=1$ to~$k.$
\end{prop}

A nice property of these linear combinations is that the coefficients
$\a\_k$ do not depend on $m$ (similarly to Proposition~\ref{teo.allk}).
For instance,
\[
L\;\,=    -(Q\_1-\si\_1I),
\]\vspace{-2.5ex}
\[
L^2=\quad\, Q\_2-\si\_1Q\_1-(\si\_2-\si_1^2)I,
\]\vspace{-2.5ex}
\[
L^3=      -(Q\_3-\si\_1Q\_2-(\si\_2-\si_1^2)Q\_1-(\si\_3-2\si\_2\si\_1+\si_1^3)I),
\]\vspace{-2.5ex}
\[
L^4=\quad\, Q\_4-\si\_1Q\_3-(\si\_2-\si_1^2)Q\_2-(\si\_3-2\si\_2\si\_1+\si_1^3)Q\_1
\]\vspace{-2.9ex}
\[
{\hspace{18em}}-(\si\_4-2\si\_3\si\_1-\si_2^2+3\si\_2\si_1^2-\si_1^4)I.
\]

\propro{\ref{pro-combi}}{
We first prove, by induction on $m,$ the identity
\beq
\label{iden-L}
(-L)^m=\suml_{k=0}^m\a'_kQ\_{m-k}
\eeq
with $\a'_0=1$ and
\beq
\label{alprim}
\a'_k=\sum_{(\b(1)\cdc\b(n\_{\b})):\;
\sum\b(i)=k}
\prod\bigl(-\si\_{\b(i)}\bigr),\quad k=1\cdc m,
\eeq
where $\b(i)$ are positive integers, $n\_{\b}$ is the variable number
of entries in $(\b(1)\cdc\b(n\_{\b})),$ and the unmarked sum and
product range from $i=1$ to ${n\_{\b}}$.

For the basis of induction, observe that $(-L)^0=I=\a'_0Q\_0.$ Let
(\ref{iden-L})--(\ref{alprim}) be valid for $(-L)^0\cdc(-L)^{m-1}$. By
Proposition~\ref{teo.allk},
\beq
\label{passa}
(-L)^{m}=\a'_0Q\_{m}-\suml_{i=0}^{m-1}\si\_{m-i}(-L)^{i}.
\eeq
Substituting (\ref{iden-L}) in the right-hand side of
(\ref{passa})
and interchanging the two sums we obtain:
\[
(-L)^{m}
=\a'_0Q\_{m}+\suml_{k=1}^m\a_{k}^{(m)}Q\_{m-k},
\]
where
\[
\a_{k}^{(m)}=\suml_{i=1}^{k}(-\si\_i)\a'_{k-i},\quad k=1\cdc m.
\]
It is easily seen that $\a_{k}^{(m)}=\a'_{k},\;$ $k=1\cdc m,$
thereby the induction step has succeeded.

Next, for an arbitrary positive integer $k,$ consider any vector
$(\b(1)\cdc\b(n\_{\b}))$ with positive integer entries such that
$\sum_{i=1}^{n\_{\b}}\b(i)=k$ (see~(\ref{alprim})).
Let $p\_j=\abs{\{i:\b(i)=j\}},$ $j=1\cdc k.$
Classifying the set of vectors $(\b(1)\cdc\b(n\_{\b}))$ such that
$\sum\b(i)=k$ by the equality of the corresponding vectors
$(p\_1\cdc p\_k),$ we see that every such a class contains
${\bigl(\sum p\_i\bigr)!}/{\prod\bigl(p\_i!\bigr)}$ members.
This implies that $\a\_k=\a'_k,\;$ $k=0,1,\ldots,$ (cf.\ (\ref{e_combi})
and~(\ref{alprim})) and thus, completes the proof.
}

\def\baselinestretch{1.19}
\section{The matrix of maximum in-forests}
\label{s_maxi}

In this section, we study some properties of the normalized matrix
$\J=J\_{n-\di}$ of maximum in-forests. Let $\mu\_{\la}(A)$ stand for the
multiplicity of $\la$ as the eigenvalue of a square matrix~$A.$

\begin{prop}
\label{LJ=0}
{\rm (i)}   ${\displaystyle L\J=\J L=LQ\_{n-\di}=Q\_{n-\di}L=0;}$\\
{\rm (ii)}  ${\displaystyle \J J\_k=J\_k\J=\J,\;\; k=\0n-\di;}$\\
{\rm (iii)} ${\displaystyle \J\mbox{\ is a projection$:$\ } \J^2=\J;}$\\
{\rm (iv)}  ${\displaystyle \rank \J=\mu\_1(\J)=\tr \J=\di;\;\;
             \mu\_0(\J)=n-\di.}$

\end{prop}

\propro{\ref{LJ=0}}{
(i) Putting $k=n-\di$ in (\ref{Fadd1}) and using the facts
that $Q\_{n-d+1}=0$ and $\si\_{n-d+1}=0,$ we get $LQ\_{n-\di}=0.$ The
other identities follow from Corollary~\ref{cor_comm} and~(\ref{Jk}).

(ii) Multiplying (\ref{j_k+1}) by $\J$ and using item~(i) and
Corollary~\ref{cor_comm}, we get the required statement, whose special
case is~(iii).

(iv) Each maximum in-forest of $\G$ has $\di$ roots, hence
$\tr Q\_{n-\di}=\di\,\si_{n-\di}$ and
$\tr \J=\tr(\si_{n-\di}^{-1}Q\_{n-\di})=\di.$ Since $\J$ is
idempotent, $\rank \J=\mu\_1(\J)=\tr \J,$ so $\mu\_0(\J)=n-\di.$
\qed
}
\smallskip

The following connection between the spectra of $L$ and
$L+\aa\J,\,$ $\aa\in\C,$ will be used in the sequel.

\begin{prop}
\label{L+J}
{\rm (i)}
The spectrum of $L+\aa\J$ consists of all nonzero eigenvalues of $L$
with their multiplicities and $\aa$ with $\mu\_{\aa}(L+\aa\J)=\di.$
{\rm (ii)}
$L+\aa\J$ is nonsingular whenever $\aa\ne0.$
\end{prop}

\propro{\ref{L+J}}{
(i)
Let $p(\la)=\si_{n-\di}^{-1}\,\sum_{i=0}^{n-\di}\si\_{n-\di-i}(-\la)^i.$
By Proposition~\ref{teo.allk}, ${\J=\si_{n-\di}^{-1}Q\_{n-\di}=p(L),}$
so $L+\aa\J=L+\aa p(L).$ Therefore, by
\cite[Theorem~3 in Chapt.~4]{Gantmacher66}, all eigenvalues of
$L+\aa\J$ are $\la'_i=\la\_i+\aa p(\la\_i),$ where $\la\_i,$ $i=\1n,$ are all
eigenvalues of $L$ with their multiplicities.
By (i)~of Proposition~\ref{LJ=0}, $L\J=0=Lp(L),$ whence $\la p(\la)$ is
an annihilating polynomial for~$L.$ Therefore, for each $\la\_i,$ a
nonzero eigenvalue of~$L,$ we have $p(\la\_i)=0,$ hence
$\la'_i=\la\_i.$ Otherwise, if $\la\_i=0,$ then
$\la'_i=\aa,$ since $p(0)=1$ by definition of~$p(\la)$. Finally, by
Corollary~\ref{mult0}, $\mu\_0(L)=\di,$ thus $\mu\_{\aa}(L+\aa\J)=\di.$
This implies~(ii).
\qed
}

\label{qnc}
\begin{prop}
\label{Tl}
$\displaystyle{\vj
=\lim_{\tau\to\infty}J(\tau)
=\lim_{\tau\to\infty} (I+\tau\,L)^{-1}.}$
\end{prop}

\propro{\ref{Tl}}{
Using Theorem~\ref{mft0}$'$ and the definition (\ref{Jtau}) of
$J(\tau),$ we have
\beqq
\lim_{\tau\to\infty} (I+\tau\,L)^{-1}
&\!\!\!=&\!\!\!
\lim_{\tau\to\infty}J(\tau)=
\lim_{\tau\to\infty}\Bigl(\suml_{k=1}^{n-d}\si\_k\tau^k\Bigr)^{-1}
                           \suml_{k=1}^{n-d}  Q\_k\tau^k
\cr
&\!\!\!=&\!\!\!
 \lim_{\tau\to\infty}\Bigl(\suml_{k=1}^{n-d}\si\_k\tau^{k-n-d}\Bigr)^{-1}
                           \suml_{k=1}^{n-d}  Q\_k\tau^{k-n-d}
=\si^{-1}_{n-d}Q\_{n-d}=\J.
\nonumber
\eeqq
}

\bigskip
\def\baselinestretch{1.15}
\section{$L$ and $\J$ as ``complementary'' linear transformations}
\label{line}

For a complex matrix $A$, let $\RR({A})$ and $\NN({A})$ denote its
range and null space, respectively. Recall that the {\em index\/} of a
square matrix $A,$ $\ind A,$ is
the smallest nonnegative integer $k$ for which $\rank(A^{k+1})=\rank(A^k).$
The {\em eigenprojection\footnote{The eigenprojections are also called {\em
principal idempotents\/} \cite{Wedderburn34,Hartwig76}.}
at\/ $0$ of $A$} \cite{Rothblum76a} or, for
short, the {\em eigenprojection of $A$} \cite{Rothblum76} is the
idempotent matrix $B$ such that $\RR(B)=\NN(A^{\nu})$ and
${\NN}(B)=\RR(A^{\nu}),$ where $\nu=\ind A.$ In other words,
$B$ is the projection {\em on $\NN(A^{\nu})$ along
$\RR(A^{\nu}).$}
The eigenprojection is
unique, because an idempotent matrix is uniquely determined by its
range and null space (see, e.g.,
\cite[p.~50]{Ben-IsraelGreville74}).\footnote{Note that for every
$A\in\C^{n\!\times\!n}$ s.t.\ $\ind A=\nu$ and every idempotent matrix
$B,$ each of the following conditions is equivalent
to $B$ being the eigenprojection of~$A$:\\
(i) $\RR(B)  =\NN(A^{\nu})$ and
    $\RR(B^*)=\NN((A^*)^{\nu})$
    \cite{Rothblum76a};\\
(ii) $A^{\nu}B=BA^{\nu}=0$ and
    $\rank A^{\nu}+\rank B=n$~\cite{Wei96,Zhang01}.\\
(iii) $AB=BA$ and $A+\aa B$ is nonsingular for all $\aa\ne0$
    \cite{KolihaStraskraba99} 
    (cf.\ (ii) of Proposition~\ref{L+J});\\
(iv) $AB=BA,$ $A+\aa B$ is nonsingular for some $\aa\ne0,$
    and $AB$ is nilpotent~\cite{KolihaStraskraba99};\\
(v) $AB=BA,$ $AB$ is nilpotent, and $AU=I-B=VA$ for some
    $U,V\in\C^{n\!\times\!n}$ \cite{Harte84};\\
(vi) $B$ commutes with all matrices commuting with $A,$
     $AB$ is nilpotent, and
     $B\ne0$ if $A$ is singular~\cite{Koliha01};\\
Moreover, the eigenprojection of $A$ is $I-AA^D,$ where $A^D$ is the
Drazin inverse of $A$ (see Section~\ref{pseudo2}).
}

Since $L\J=0$ (Proposition~\ref{LJ=0}), we have $\RLT\cap\RJ=\{{\bf
0}\},$ where $L^*=L^{\interca}.$
Similarly, $\J L=0$ implies $\RJT\cap\RL=\{{\bf
0}\}.$ Consequently, by \cite[Theorem~11]{MarsagliaStyan74}, $L$ and
$\J^*$ are {\em rank additive}, i.e., $\rank(L+\J^*)=\rank L+\rank \J^*.$
Corollary~\ref{mult0} implies that $\rank L\ge n-\di,$ whereas, by
Proposition~\ref{LJ=0}, $\rank \J^*=\di.$ Since $\rank(L+\J^*)\le n,$
we have $\rank L=n-\di$ and $\rank(L+\J^*)=n.$
Now $L\J=\J L=0$ implies $\NL=\RJ$ and $\NJ=\RL.$
Furthermore, by Proposition~\ref{L+J}, $\rank(L+\J)=n,$ hence
$L$ and $\J$ are rank additive. It follows now from
\cite[Theorem~11]{MarsagliaStyan74} that $\RL\cap\RJ=\{{\bf 0}\}.$
Since $\RJ=\NL,$ we get $\RL\cap\NL=\{{\bf 0}\},$ which, by
\cite[p.~165]{Ben-IsraelGreville74}, implies $\ind L=1.$
The latter fact together with $\RJ=\NL,$ $\NJ=\RL,$ and $\J^2=\J$ imply
that $\J$ is the eigenprojection of~$L$ (alternatively, this follows
from Proposition~\ref{Tl} and \cite[Theorem~3.1]{Meyer74}).
We proved
\begin{prop}
\label{L+JT}
{\rm(i)} $L+\J^*$ is nonsingular.\\
{\rm(ii)} $\;\rank L=n-\rank \J=n-\di.$\\
{\rm(iii)} $\NL=\RJ$ and\/ $\RL=\NJ.$\\
{\rm(iv)} $\RL\cap\RJ=\{{\bf 0}\}.$\\
{\rm(v)} $\;\ind L=1.$\\
{\rm(vi)} $\J$ is the eigenprojection of~$L.$
\end{prop}

It is known \cite[p.~194]{Rothblum76a}, \cite[Theorem~7.a.3]{Rothblum81SIAM}
that for every finite homogeneous Markov chain with a transition
matrix~$P,$ the {\em long run transition matrix}
$\Pin=\liml_{k\to\infty}\frac{1}{k}\,\suml_{t=0}^{k-1} P^t$ is the
eigenprojection of $P$ at $1,$ which is
the eigenprojection of $I-P.$\footnote{This also follows from Meyer's
Theorem~2.2 in~\cite{Meyer75}. Indeed, by this theorem,
$\Pin\!=\!I-(I-P)(I-P)^{\gri},$ where $(I-P)^{\gri}$ is the group
inverse of $I-P,$ and the right-hand side is the eigenprojection of
$I-P,$ as mentioned in the next section.} On the other hand, $I-P$ is
exactly the Laplacian matrix $L$ of the weighted digraph without loops
whose arc weights are equal to the corresponding transition
probabilities. Therefore $\J,$ the eigenprojection of $L,$ coincides
with~$\Pin.$ The fact that $\Pin$
coincides with the normalized matrix of maximum in-forests of the
digraph corresponding to a Markov chain is the so called {\em
Markov chain tree theorem\/} \cite{LeightonRivest83,LeightonRivest86}.
Thus, item (vi) of Proposition~\ref{L+JT} provides an immediate proof
of this theorem.

By virtue of Proposition~\ref{LJ=0}, every nonzero column of $\J$ (or
$Q\_{n-\di}$) is an eigenvector of $L$ that corresponds to the zero
eigenvalue. Moreover, it follows from $\NL=\RJ$
(Proposition~\ref{L+JT}) that the nonzero columns of $\J$ span the null
space of~$L$.
Since, by (\ref{adjexp}), $Q(\tau)$ is proportional to $\adj(\la I-(-L))$
at $\la=\tau^{-1},$ $Q(\tau)$ can be used to generate some eigenvectors
of $L$ that correspond to its nonzero eigenvalues. For completeness, we
give a proof of this fact.

\begin{prop}
\label{evect}
Let $\la\_i\ne0$ be an eigenvalue of~$L.$ Then every nonzero column of
$Q(-\la_i^{-1})$ is an eigenvector of $L$ that corresponds to~$\la\_i.$
\end{prop}

\propro{\ref{evect}}{
Let $X=\la\_i\!I-L.$ Then $\det X=0$. Using Theorem~\ref{mft0}$'$ and
the fact that for every square matrix $Y,\;$ $Y\adj Y=(\det Y)I$
holds, we get
\[
(\la\_i\!I-L)\,Q(-\la_i^{-1})
= X\adj(I-\la_i^{-1}L)
=\la_i^{1-n}X\adj X
=\la_i^{1-n}(\det X)I=0.
\]
This implies the desired statement.
}

\section{Forest matrices and generalized inverses of~$L$}
\label{pseudo2}

The {\it Moore-Penrose generalized inverse\/} $A^+$ of a
rectangular complex matrix $A$ is the unique matrix $X$ such that
\[
AXA=A,\quad XAX=X,\quad (AX)^*=AX,\quad (XA)^*=XA.
\]

For an arbitrary square matrix $A,$ its {\em Drazin inverse},
$A^D,$ is the unique matrix $X$ satisfying the equations
\[
A^{\nu+1}X=A^{\nu},\quad XAX=X,\quad AX=XA,
\]
where $\nu=\ind A.$
If $\nu=0,$ then $A^D=A^{-1}$; if $\nu\le1,$ then $A^D$ is referred to
as the {\em group inverse}, $A^{\gri},$ i.e., the unique matrix $X$
such that
\[
AXA=A,\quad XAX=X,\quad AX=XA.
\]

As applied to the Laplacian matrices of graphs, the generalized
inverses were considered in connection with the analysis of electrical
networks (providing ``resistance distance''), Markov chains, and some
preference aggregation problems (more specifically, estimation from
paired comparisons), in constructing geometrical representations of graphs
(with applications to chemistry, social networks, etc.), in control,
cluster analysis, and parallel computing.
There is a huge literature on generalized inverses within the last
years. For multiple representations of
the Drazin inverse, see, e.g.,~\cite{WeiWu00,ChenChen00,Chen01}.

In this section, we present a few relations between the $\Lpr$ and
the forest matrices and one representation for $L^+$. In the case
of symmetric $L,$ where $\Lpr=L^+,$ some of these expressions are
given in~\cite{CheSha98}.\footnote{For symmetric $L,$ interesting
representations for $\Lpr=L^+$ were proposed in \cite{Fiedler95},
\cite[Theorem~2.2]{KirklandNeumann97}, and, in case of weighted
trees, in~\cite{KirklandNeumann97} and~\cite[Theorem~3]{Bapat97}.
In~\cite[Theorem~3]{DruryStyan94} a combinatorial interpretation
of the Campbell-Youla inverse (the symmetric generalized inverse
with the zero diagonal) of~$L$ is given.}

For an arbitrary square matrix $A,\,$ $AA^D$ is the unique projection
on $\RR(A^{\nu})$ along $\NN(A^{\nu})$
\cite[p.~173]{Ben-IsraelGreville74}. Then $I-AA^D$ is the projection
on $\NN(A^{\nu})$ along $\RR(A^{\nu}).$
Therefore, $I-AA^D$ is the eigenprojection of
$A$~\cite{Rothblum76a,Rothblum76}. Combining this with items (v) and
(iv) of Proposition~\ref{L+JT}, we obtain

\begin{prop}
\label{p_eig-proj}
$\J=I-LL^{\gri}.$
\end{prop}

The fact that $\J$ is the eigenprojection of $L$ helps interpret, in
terms of in-forests, the expressions of generalized inverses of $L$
that involve the eigenprojection of~$L.$

\begin{prop}
\label{Lgroupinv}
$\mathstrut$\\
{\rm (i)}$\;$ For any $\aa\ne 0,\;$ $\Lpr=(L+\aa\vj)^{-1}-\aa^{-1}\vj,$
whence $\Lpr=\liml_{\abs{\aa}\to\infty}(L+\aa\vj)^{-1}.$\\
{\rm (ii)}    For any $\aa\ne 0,\;$ $\Lpr=(L+\aa\vj)^{-1}(I-\vj)$.\\
{\rm (iii)} $\displaystyle{\Lpr=\frac{\mathstrut\si_{n-\di-1}}{\si_{n-\di}}
            \left(J_{n-\di-1}-\q\right).}$\\
{\rm (iv)}  $\displaystyle{\Lpr
            =\liml_{\tau\to\infty}\tau\left(J(\tau)-\vj\right).}$
\end{prop}

{\bf Remarks on Proposition~\ref{Lgroupinv}.}
(i), (iii), and (iv) were presented in~\cite{AgaChe01}.
(i)~results by substituting $\J$ for the eigenprojection in the
expression of group inverse employed
in~\cite[p.~150]{MeyerStadelmaier78} (for its proof
see~\cite[Theorem~4.2]{Rothblum81LAA}; related expressions appeared
in \cite[Theorem~5.5]{Meyer75} and \cite[last line on
p.~646]{Rothblum76}, where `$+$' must be replaced by `$-$').
(ii)~is obtained by the same substitution in the representation of
Drazin inverse given in \cite{Koliha01} (the case with $\aa=1$ appeared
in~\cite{Rothblum76}) or by multiplying (i) by $LL^{\gri}=I-\J.$
In view of Propositions~\ref{char-L} and~\ref{MatCo}, (iii) follows
from the expression of Drazin inverse discovered independently by
Hartwig~\cite[Eq.~(13)]{Hartwig76} and Gower~\cite[Theorem~1]{Gower80}.
\medskip

The matrices $L+\aa\J$ are the ``complementary
perturbations''~\cite{MeyerStadelmaier78} of~$L.$ Matrices of this kind
are important for the analysis of M-matrices and singular systems of
equations. In particular, a matrix $A$ with eigenprojection $B$ and
nonpositive off-diagonal entries is an M-matrix if and only if for some
$c>0,\,$ $(A+\aa B)^{-1}$ is nonnegative when
$\aa\in(0,c)$~\cite{MeyerStadelmaier78}. If $A$ is an M-matrix, then
$(A+\aa B)^{-1},$ $\aa\in(0,c),$ make up a class of {\em nonnegative
nonsingular commuting weak inverses\/} for
$A$~\cite{MeyerStadelmaier78}.  $(L+\aa\J)^{-1}$ can be represented as
a linear combination of forest matrices using (i) and (iii) of
Proposition~\ref{Lgroupinv}:
\[
(L+\aa\J)^{-1}
=\frac{\si_{n-\di-1}}{\si_{n-\di}}\left(J_{n-\di-1}+\b\J\right),
\]
where $\b=\frac{\si_{n-\di}}{\aa\si_{n-\di-1}}-1.$ This throws some
light on the nonnegativity of $(L+\aa\J)^{-1}$: if
$\aa\in(0,\frac{\si_{n-\di}}{\si_{n-\di-1}})$ then $(L+\aa\J)^{-1}$ is
a positive combination of $J_{n-\di-1}$ and $\J.$ Based on this, we
termed $(L+\aa\J)^{-1}$ the {\em matrices of dense in-forests\/}
of~$\G.$ These and the inverse ``uniform diagonal perturbations''
$(L+\aa I)^{-1}$ can serve to measure proximity between digraph
vertices~\cite{AgaChe01}. Note in this connection that by
\cite[Corollary~4.4]{Rothblum81LAA}, $(L+\aa\J)^{-1}_{ij}>0$ for all
$\aa>0$ sufficiently small if and only if vertex $j$ is accessible from
$i$ in $\G,$ and the same is true for $(L+\aa I)^{-1}_{ij}.$ By
Theorem~\ref{mft0}$',$ $(L+\aa I)^{-1}$ is proportional to $J(\tau)$
with $\tau=1/\aa.$

\medskip

We conclude with one expression for the Moore-Penrose inverse of~$L.$

Consider the matrix $Z\!:=L+\J^{*}$ which is nonsingular by
Proposition~\ref{L+JT}. Using the identity $L\J=0$
(Proposition~\ref{LJ=0}), we get
$
{(Z^{*})}^{-1}Z^{-1} =(ZZ^{*})^{-1}
=(\J^{*}\!\!\vj+LL^{*})^{-1}.
$

\def\afterthmseparator{}
\begin{prop}
\label{pseudoor}
{\rm \cite{AgaChe01}.}
$L^+
=L^{*}{(ZZ^{*})}^{-1}
=L^{*}(\J^{*}\!\!\vj+L L^{*})^{-1}.$
\end{prop}
\def\afterthmseparator{.}

One method to prove this is to check the conditions in
the definition of Moore-Penrose inverse by direct computation using
Proposition~\ref{LJ=0} and the
facts that ${(ZZ^{*})}^{-1}$ commutes with $LL^{*}$ and
$\J^{*}\!\!\vj$ and that
$LL^{*}{(ZZ^{*})}^{-1}$ and
$\J^{*}\!\!\vj{(ZZ^{*})}^{-1}$ are symmetric~\cite{AgaChe01}.
Alternatively, Proposition~\ref{pseudoor} can be proved by
employing the Penrose formula $A^+=A^*(AA^*)^+$, the fact that
$(AA^*)^+=(AA^*)^{\gri}$ (since $AA^*$ is Hermitian) and an expression
of $(AA^*)^{\gri}$ such as those given in (i) and (ii) of
Proposition~\ref{Lgroupinv}.

\section{A concluding remark}

It is instructive to compare the ``Laplacian graph mathematics'' we
touched upon in this paper with the corresponding results on the
adjacency characteristic matrix, see, e.g., \cite[Sections 1.4, 1.9.1,
1.9.5 and others]{CvetkovicDoobSachs80} and the articles by Kasteleyn
and Ponstein cited therein, \cite{Schwenk91}, and so on. This
comparison suggests that the Laplacian mathematics is based on trees in
the same sense as the ``adjacency graph mathematics'' is based on
routes and circuits. We mean that a number of expressions related with
the adjacency characteristic matrix can be interpreted in terms of
routes and circuits, whereas the counterparts of these expressions
related with the Laplacian characteristic matrix involve spanning
forests for their interpretation.
\vspace{-.9em}

\section*{Acknowledgements}

This work was supported in part by the Russian Foundation for Basic
Research under grants 01--01--10732-z and 02--01--00614-a. The authors
thank a referee for a careful reading of the paper and many thoughtful
suggestions.

{\small

}
\end{document}